\documentstyle[amsfonts,11pt]{article}
\newcommand{\Y}{\mbox{$Y$}}
\newcommand{\M}{\mbox{$M$}}

\newcommand{\OM}{\mbox{$\Omega^{1}M$}}

\newcommand{\X}{\mbox{$X$}}
\newcommand{\p}{\partial}
\newcommand{\f}{\cal O}
\newcommand{\ko}{{\sf k}}
\renewcommand{\i}{\mathrm i}
\renewcommand{\j}{\mathrm j}
\newcommand{\k}{\mathrm  k}
\newcommand{\lon}{\mbox{$\longrightarrow$}}
\newcommand{\hook}{\mbox{$\hookrightarrow$}}

\newcommand{\bnu}{\nu}
\newcommand{\bmu}{\mu}

\newtheorem{theorem}{Theorem}

\newtheorem{corollary}[theorem]{Corollary}

\begin{document}

\title{Projective structures on moduli spaces of compact complex
hypersurfaces }

\author{Sergey Merkulov \\
\small School of Mathematics and Statistics, University of Plymouth\\
\small Plymouth, Devon PL4 8AA, United Kingdom\\
       Henrik Pedersen\\
\small   Department of Mathematics and Computer Science, Odense University, \\
\small                Campusvej 55, 5230 Odense M, Denmark
}

\date{}
\maketitle

\begin{abstract}
It is shown that moduli spaces of complete families of compact complex
hypersurfaces in complex manifolds often come equipped canonically with
projective structures satisfying some natural integrability conditions.
\end{abstract}


\sloppy
\paragraph{1. Projective connections.}
Let $M$ be a complex manifold. Consider the following equivalence
relation on the set of affine torsion-free connections on $M$: two connections
$\hat\Gamma$ and $\Gamma$    are said to be projectively equivalent
if they have the same geodesics, considered as unparameterized
paths. In a local coordinate chart $\{t^{\alpha}\}$, $\alpha =
1,\ldots,\dim M$ on $M$, where $\hat{\Gamma}$ and $\Gamma$ are represented by
Christoffel symbols $\hat{\Gamma}_{\alpha\beta}^{\gamma}$ and
${\Gamma}_{\alpha\beta}^{\gamma}$ respectively,
 this equivalence relation reads \cite{Hitchin}
$$
\hat{\Gamma} \sim {\Gamma}
$$
if
$$
\hat{\Gamma}_{\alpha\beta}^{\gamma} = \Gamma_{\alpha\beta}^{\gamma}
+ b_{\beta}\,\delta_{\alpha}^{\gamma} + b_{\alpha}\,\delta^{\gamma}_{\beta}
$$
for some 1-form $b = b_{\alpha}\,dt^{\alpha}$. An equivalence class of
torsion-free affine connections under this relation is called a
projective structure or a projective connection.

Let $M$ be a complex manifold with a projective structure.
A complex submanifold $P\subset M$ is called
{\em totally geodesic}
if for each point $t\in P$ and each direction tangent to $P$ at $t$ the
corresponding geodesic of the projective connection is contained in $P$ at
least locally.

\paragraph{2. Moduli spaces of compact complex hypersurfaces.}
Let $X$ be a compact complex hypersurface in a complex manifold $Y$ with
normal line bundle~$N$ such that
$H^{1}(\X,N) = 0$. According to Kodaira \cite{Kodaira1},  such a
 hypersurface \X\ belongs to the complete analytic
 family \{$\X_{t} \mid t \in \M$
\} of compact complex hypersurfaces $\X_{t}$ in $Y$ with the moduli space $M$
being
a $\left( \dim_{\Bbb C}H^{0}(\X, N)\right)$-dimensional complex manifold.
Moreover, there is a canonical isomorphism
$\ko_{t}: T_{t}\M \longrightarrow  H^{0}(X_{t}, N_{t})$
which associates a global
section of the normal bundle $N_{t}$ of $\X_{t} \hookrightarrow Y$ to any
tangent vector at the corresponding point $t \in \M$.

 Consider
$F= \{(y,t)\in Y\times M \mid y\in X_t\}$ and denote
by $\bmu: F\lon Y$ and $\bnu: F\lon M$ two natural projections,
\begin{equation}
Y \stackrel{\bmu}{\longleftarrow} F \stackrel{\bnu}{\longrightarrow} M .
\label{F1}
\end{equation}
The space $F$ is
 a submanifold in  $Y\times M$.
If $N_{F}$ is the normal bundle of $F \hookrightarrow \Y\times \M$, then,
for any point $t \in \M$, we have
$\left. N_{F}\right|_{\bnu^{-1}(t)}\, \simeq N_{X_{t}\mid Y}$,
where $N_{X_{t}\mid Y}$ is the normal bundle of the submanifold
$\bmu\circ\bnu^{-1}(t) =
X_{t} \hookrightarrow Y$. By Kodaira's theorem, there is an isomorphism
$\ko: TM {\longrightarrow} \bnu_{*}^{0}(N_{F})$,
where $\bnu_{*}^{0}(N_{F})$ denotes the direct image.

Let us denote the point in the
moduli space $M$ corresponding to $X$ by $t_{0}$, i.e.\
$X=\bmu\circ\bnu^{-1}(t_{0})$. It is easy to show that,
for each $y \in Y' \equiv \cup_{t\in M}
X_{t}$  the set $\bnu\circ\bmu^{-1}(y)$ is a complex analytic
subspace of $M$.  We denote by $P_{y}$ its manifold content,
i.e. $P_{y} = \bnu\circ\bmu^{-1}(y)\setminus\{${\small singular points}$\}$.
If the natural evaluation map
\begin{eqnarray*}
H^{0}(X_{t}, N_{X_{t}\mid Y}) & \lon & N_{z} \\
\phi &\lon & \phi(z) ,
\end{eqnarray*}
where $N_{z}$ is the fibre of  $N$ at a point $z\in X_{t}$
and $\phi(z)$ is the value of the  global section $\phi\in
H^{0}(X_{t}, N_{X_{t}\mid Y})$  at $z$, is surjective
at all points $z\in X_{t}$ and for all $t\in M$, then
$P_{y} = \bnu\circ\bmu^{-1}(y)$.

\paragraph{3. The main theorem.}
The idea of studying differential geometry on the
moduli space of compact complex submanifolds of a given ambient complex
manifold goes back to Penrose \cite{Penrose} who discovered self-dual conformal
structures automatically induced on 4-dimensional moduli spaces of rational
curves with normal bundle $N={\Bbb C}^2\otimes \f (1)$.  In this section we
show that moduli spaces of compact complex hypersurfaces often come equipped
canonically with induced projective structures satisfying some natural
integrability conditions.  Other manifestations of general and strong links
between complex analysis and differential geometry can be found in Merkulov's
survey \cite{Merkulov}.

\begin{theorem}\label{th 1}
Let \mbox{$X \hookrightarrow Y$} be a compact complex submanifold of
codimension 1 with normal
bundle $N$ such that $H^{1}(X, N) = 0$ and let $M$ be the associated
complete moduli space of
relative deformations of $X$ inside $Y$. If $H^{1}(X,\f_{X}) = 0$, then
a sufficiently small neighbourhood $M_{0}\subset M$ of the
point $t_{0}\in M$ corresponding to $X$, comes
equipped canonically with a projective structure
such that, for every point $y \in Y' \equiv \cup_{t\in M}
X_{t}$, the associated submanifold $P_{y}\subseteq \bnu\circ \bmu^{-1}(y) \cap
M_{0}$ is totally geodesic.
\end{theorem}
{\em Proof}.  An open neighbourhood of the submanifold
$X\hookrightarrow Y$ can
always be covered by  a finite number of coordinate
charts $\{W_{\i}\}$ with local coordinate functions $(w_{\i},
z^{a}_{\i})$, $a=1,\ldots,n=\dim X$, on each
neighbourhood $W_{\i}$ such that $X\cap W_{\i}$ coincides with the subspace of
$W_{\i}$ determined by the equation $w_{\i} = 0$. On the intersection
$W_{\i} \cap W_{\j}$
the coordinates $w_{\i}, z^{a}_{\i}$ are holomorphic functions of
$w_{\j}$ and $z^{b}_{\j}$,
$$
w_{\i} = f_{\i\j}(w_{\j}, z^{b}_{\j}),\ \ \ \
z^{a}_{\i} = g^{a}_{\i\j}(w_{\j}, z^{b}_{\j}),
$$
with $f_{\i\j}(0, z^{b}_{\j}) = 0$. Here $z_{\j}=(z_{j}^{1},\ldots,
z^{n}_{j})$.

Let $U\subset M$ be a coordinate
neighbourhood of the point $t_{0}$ with
coordinate functions
$t^{\alpha}, \alpha = 1,\ldots,m=\dim M $. Then the
coordinate domains
$U\times W_{\i}$ with coordinate functions $\left(w_{\i}, z^{a}_{\i},
t^{\alpha}\right)$ cover an open neighbourhood of $X\times U$ in the manifold
$Y\times U$. For a sufficiently small $U$,
the submanifold $F_{U} \equiv \bnu^{-1}(U)\hookrightarrow Y\times U$ is
described in each
coordinate chart $W_{\i}\times U$ by an equation of the form \cite{Kodaira1}
$$
w_{\i} = \phi_{\i}(z_{\i}, t),
$$
where  $\phi_{\i}(z_{\i}, t)$ is a holomorphic function
of $z^{a}_{\i}$ and $t^{\alpha}$ which satisfies the boundary conditions
$\phi_{\i}(z_{\i}, t_{0}) = 0$.
For each fixed $t\in U$ this equation defines a submanifold
$X_{t}\cap W_{\i}\hookrightarrow W_{\i}$.

By construction, $F_{U}$ is covered by a finite number of
coordinate neighbourhoods
$\{V_{\i}\equiv \left.W_{\i}\times U\right|_{F}\}$ with local coordinate
functions $(z^{a}_{\i},
t^{\alpha})$
which are related to each other on the intersections $V_{\i}\cap V_{\j}$
as follows
$$
z^{a}_{\i} = g^{a}_{\i\j}\left(\phi_{\j}(z_{\j},t), z_{\j}\right).
$$
Obviously we have
$\phi_{\i}\left(g_{\i\j}\left(\phi_{\j}(z_{\j},t),
z_{\j}\right), t\right) = f_{\i\j}\left(\phi_{\j}(z_{\j}, t), z_{\j}\right)$.

The Kodaira map $\ko: \left.TM\right|_{U} \lon
\bnu_{*}^{0}(\left.N_{F}\right|_{F_{U}})$ can be described in the
following way: take
any vector field $v$ on $U$ and apply the corresponding
1st-order differential operator $ V^{\alpha}\p_{\alpha}$, where $\p_{\alpha} =
\p/ \p t^{\alpha}$, to each function
$\phi_{\i}(z_{\i}, t)$. The result is a collection
of holomorphic functions
$\sigma_{\i}(z_{\i},t) =  V^{\alpha}\,
\partial_{\alpha} \phi_{\i}(z_{\i}, t)$
defined respectively on $V_{\i}$. On the intersection $V_{\i}\cap V_{\j}$
one has $\left.\sigma_{\i}(z_{\i},t)\right|_{z_{\i}= g_{\i\j}(\phi_{j}, z_{j})}
= F_{\i\j}(z_{\j},t)\,
\sigma_{j}(z_{\j},t)$,
where
$$
F_{\i\j} \equiv
  \left.\frac{\partial f_{\i\j}}{\partial w_{\j}}\right|_
{w_{\j} = \phi_{\j}(z_{\j},t)}  -
\left.\frac{\p \phi_{\i}}{\p z_{\i}^{a}}\right|_{z_{\i} =
g_{\i\j}(\phi_{\j}, z_{\j})}
\left.\frac{\p g_{\i\j}^{a}}{\p w_{\j}}\right|_{w_{\j} =
\phi_{\j}(z_{\j},t)},
$$
is the transition matrix of the normal bundle $\left.N_{F}\right|_{F_{U}}$
on the  overlap $F_{U}\cap V_{\i}\cap V_{\j}$. Therefore the 0-cochain
$\left\{\sigma_{\i}(z_{\i}, t) \right\}$
is  a \v{C}ech 0-cocycle
representing a global section $\ko(v)$ of the normal bundle $N_{F}$ over
$F_{U}$.

Let us investigate how second partial derivatives of $\{\phi_{\i}(z_{\i},t)\}$
and $\{\phi_{\j}(z_{\j},t)\}$
are related on the intersection $V_{\i}\cap V_{\j}$. Since
$$
\left.\frac{\partial\phi_{\i}(z_{\i},t)}{\partial t^{\alpha}}\right|_{z_{\i}
= g_{\i\j}(\phi_{\j}, z_{\j})} = F_{\i\j}\,
\frac{\partial\phi_{\j}(z_{\j},t)}{\partial t^{\alpha}}
$$
we find
\begin{equation}
\left.\frac{\partial^{2}\phi_{\i}}{\partial
t^{\alpha} \partial t^{\beta}}\right|_{z_{\i}
= g_{\i\j}(\phi_{\j}, z_{\j})} =
F_{\i\j}\,
\frac{\partial^{2}\phi_{\j}}{\partial
t^{\alpha} \partial t^{\beta}}
+ E_{\i\j}\,
\frac{\partial\phi_{\j}}{\partial t^{\alpha}}
\frac{\partial\phi_{\j}}{\partial t^{\beta}}
- G_{\i\j\,\alpha}\, \frac{\partial\phi_{\j}}{\partial t^{\beta}}
- G_{\i\j\,\beta}\, \frac{\partial\phi_{\j}}{\partial t^{\alpha}},
\label{04}
\end{equation}
where
\begin{eqnarray*}
E_{\i\j} & = & \left.\frac{\partial^{2}f_{\i\j}}{\partial
w_{\j} \partial w_{\j}}\right|_{w_{\j} = \phi(z_{\j},t)} -
\left.\frac{\p \phi_{\i}}{\p z_{\i}^{a}}\right|_{z_{\i} = g_{\i\j}(\phi_{\j},
z_{\j})}
\left.\frac{\p ^{2} g_{\i\j}^{a}}{\p w_{\j}\p w_{\j}}\right|_{w_{\j} =
\phi_{\j}(z_{\j},t)}\\
&& - \, \left.\frac{\p ^{2} \phi_{\i}}{\p z_{\i}^{a}\p
z_{\i}^{b}}\right|_{z_{\i} = g_{\i\j}(\phi_{\j}, z_{\j})}
\left. \left(
\frac{\p g_{\i\j}^{a}}{\p w_{\j}}
\frac{\p g_{\i\j}^{b}}{\p w_{\j}}
\right)\right|_{w_{\j} = \phi_{\j}(z_{\j},t)},
\end{eqnarray*}
and
$$
G_{\i\j\,\alpha} =
\left.\frac{\p ^{2} \phi_{\i}}{\p z_{\i}^{a}\p t^{\alpha}}\right|_{z_{\i} =
g_{\i\j}(\phi_{\j}, z_{\j})}
\left.\frac{\p g_{\i\j}^{a}}{\p w_{\j}}\right|_{w_{\j} =
\phi_{\j}(z_{\j},t)}.
$$
The collections $\left\{E_{\i\j}\right\}$ and
$\left\{G_{\i\j\,\alpha}\right\}$ form 1-cochains with coefficients
in $N_{F}^{*}$ and $\bnu^{*}(\OM)$, respectively. Straightforward
calculations reveal the obstructions for these two 1-cochains to be 1-cocycles,
\begin{eqnarray*}
\delta \left\{E_{\i\k}\right\} & = & 2\,
\frac{\p F_{\i\j}(z_{\j},t)}{\p z_{\j}^{a}}\,
\left.\frac{\p g_{\j\k}^{a}}{\p w_{\k}}\right|_{w_{k}=\phi_{\k}(z_{\k},t)}  \\
\delta \left\{G_{\i\k\,\alpha}\right\} & = &
 \frac{\p F_{\i\j}(z_{\j},t)}{\p z_{\j}^{a}}\,
\left.\frac{\p g_{\j\k}^{a}}{\p w_{\k}}\right|_{w_{k}=\phi_{\k}(z_{\k},t)}\,
\frac{\partial\phi_{j}(z_{\j},t)}{\partial t^{\alpha}}.
\end{eqnarray*}
 From these equations we conclude that
the 1-cochain
$\left\{\tau_{\i\k\,\alpha}\right\}$, where
$$
\tau_{\i\k\,\alpha}\equiv
\frac{1}{2}E_{\i\k}\,\frac{\partial\phi_{\k}}{\partial
t^{\alpha}} - G_{\i\k\,\alpha},
$$
is actually a 1-cocycle with values in
$\bnu^{*}(\OM)$. Since $H^{1}(X, \f_{X}) = 0$, the
semi-continuity principle \cite{Kodaira2} implies $H^{1}(X_{t}, \f_{X_{t}})
= 0$ for all points in some Stein neighbourhood $M_{0}\subseteq U$. Hence, by
the Leray spectral sequence
$H^{1}\left(\bnu^{-1}(M_{0}),\bnu^{*}(\OM)\right) = 0$.
Therefore, the 1-cocycle $\left\{\tau_{\i\k\,\alpha}\right\}$ is always a
coboundary
$ \left\{\tau_{\i\j\,\alpha}\right\} =
\delta\left\{\theta_{\i\,\alpha}\right\}$,
or more explicitly,
\begin{equation}
\tau_{\i\j\,\alpha}(z_{\j},t) = F_{\i\j}(z_{j},t)\left(- \left.
\theta_{\i\,\alpha}(z_{\i},t)\right|_{z_{\i}=g_{\i\j}(\phi_{\j}, z_{\j})} +
\theta_{\j\,\alpha}(z_{\j},t)
\right),\label{05}
\end{equation}
for some 0-cochain $\left\{\theta_{\i\,\alpha}(z_{\i},t)\right\}$
on $\bnu^{-1}(M_{0})$
with values  in $\bnu^{\ast }(\OM)$. However, this
0-cochain is defined non-uniquely ---
for any global section $\xi = \xi_{\alpha}\,dt^{\alpha}$ of
$\bnu^{*}(\OM)$ over $\bnu^{-1}(M_{0})$
the 0-cochain
\begin{equation}
\tilde{\theta}_{\i\,\alpha}(z_{\i},t) = \theta_{i\,\alpha }(z_{\i},t) +
\left.\xi_{\alpha}(t)\right|_{\bnu^{-1}(M_{0})\cap V_{i}} \label{06}
\end{equation}
splits the same 1-cocycle $\left\{\tau_{\i\j\,\alpha}\right\}$. Note that,
due to the compactness of the complex submanifolds $\bnu^{-1}(t)\subset F$
for all $t\in M_{0}$ the components $\xi_{\alpha}$ of the global section
$\xi \in H^{0}(\bnu^{-1}(M_{0}), \bnu^{*}(\OM))$ are constant
along the fibers, i.e. $\xi_{\alpha}\in \bnu^{-1}( \f_{M_{0}})$.

If we rewrite equation (\ref{04}) in the form
\begin{eqnarray*}
\left.\frac{\partial^{2}\phi_{\i}(z_{\i},t)}{\partial
t^{\alpha} \partial t^{\beta}}\right|_{z_{i}
= g_{\i\j}(\phi_{\j}, z_{\j})} & = &
F_{\i\j}(z_{\j},t)\,
\frac{\partial^{2}\phi_{\j}(z_{\j},t) }{\partial
t^{\alpha} \partial t^{\beta}} \\
 && +\, \tau_{\i\j\,\alpha}(z_{\j},t)\,
\frac{\partial\phi_{\j}(z_{\j},t)}{\partial t^{\beta}}\,
+ \tau_{\i\j\,\beta}(z_{\j},t)\,
\frac{\partial\phi_{\j}(z_{\j},t)}{\partial t^{\alpha}}
 \end{eqnarray*}
and take equation (\ref{05}) into account, we obtain the
equality
$$
\left.\left(\frac{\partial^{2}\phi_{\i}}{\partial
t^{\alpha} \partial t^{\beta}}
+
\theta_{\i\,\alpha}\,\frac{\partial\phi_{\i}}{\partial t^{\beta}}
+
\theta_{\i\,\beta}\,\frac{\partial\phi_{\i}}{\partial t^{\alpha}}
\right)\right|_{z_{\i} = g_{\i\j}(\phi_{\j}, z_{\j})} =
\frac{\partial^{2}\phi_{\j}}{\partial
t^{\alpha} \partial t^{\beta}} +
\theta_{\j\,\alpha}\,\frac{\partial\phi_{\j}}{\partial t^{\beta}}
+
\theta_{\j\,\beta}\,\frac{\partial\phi_{\j}}{\partial t^{\alpha}}
$$
which implies that, for each value of $\alpha$ and
$\beta$, the holomorphic functions,
$$
\Phi_{\i\,\alpha\beta}(z_{\i},t) \equiv
\frac{\partial^{2}\phi_{\i}(z_{\i},t)}{\partial
t^{\alpha} \partial t^{\beta}}
+
\theta_{\i\,\alpha }(z_{\i},t)\,\frac{\partial\phi_{\i}(z_{\i},t)}{\partial
t^{\beta}}
+
\theta_{\i\,\beta}(z_{\i},t)\,\frac{\partial\phi_{\i}(z_{\i},t)}{\partial
t^{\alpha}},
$$
represent a global section of the normal bundle $N_{F}$ over
$\bnu^{-1}(M_{0})$. Since the collections of functions
$\left\{\partial_{\alpha}\phi_{\i}(z_{\i},t)\right\}$
form a \v{C}ech representation of a basis for the free $\f_{M_{0}}$-module
$\bnu_{*}^{0}\left(\left.N_{F}\right|_{\bnu^{-1}(M_{0})}\right)$,
the equality
\begin{equation}
\Phi_{\i\alpha\beta}(z_{\i},t) = \Gamma_{\alpha\beta}^{\gamma}(t)\,
\partial_{\alpha}\phi_{\i}(z_{\i},t) \label{07}
\end{equation}
must hold for some global holomorphic functions
$\Gamma_{\alpha\beta}^{\gamma}$ on $\bnu^{-1}(M_{0})$. Since all the fibers
$\bnu^{-1}(t)$, $t\in M_{0}$, are compact complex manifolds, these functions
are actually
pull-backs of some holomorphic functions on $M_{0}$. A coordinate system
$\{t^{\alpha}\}$ on $M_{0}$ was used in the construction of
$\Gamma_{\alpha\beta}^{\gamma}(t)$. However from (\ref{07}) it
immediately follows that under general
coordinate transformations $t^{\alpha} \longrightarrow t^{\alpha'} =
t^{\alpha'}(t^{\beta})$
these functions transform according to
$$
\Gamma_{\alpha' \beta'}^{\gamma'} = \frac{\partial t^{\gamma'}}
{\partial t^{\delta}}\left( \Gamma_{\mu\nu}^{\delta}\, \frac{\partial
t^{\mu}}{\partial t^{\alpha'}}\,
\frac{\partial t^{\nu}}{\partial t^{\beta'}} + \frac{\partial^{2}
t^{\delta}}{\partial t^{\alpha'}
\partial t^{\beta'}}\right).
$$
Thus from any given splitting $\left\{\tau_{\i\j\,\alpha}\right\} =
\delta\left\{\theta_{\i\,\alpha}\right\}$
 of the 1-cocycle $\left\{\tau_{\i\j\,\alpha}\right\}$ we extract a
symmetric affine
connection $\Gamma_{\alpha\beta}^{\gamma}(t)$.
It is straightforward to check that this connection is independent of the
choice of the $(w_{\i},z_{\i}^{a})$-coordinate system used in the construction
and thus is well-defined except for the arbitrariness in its construction
described by the transformations
(\ref{06}) which, as one can easily check, change the connection as follows
\begin{eqnarray*}
\theta_{\i\,\alpha}(z_{\i},t) & \lon & \theta_{\i\,\alpha}(z_{\i},t) +
\xi_{\alpha}(t) \\
\Gamma_{\alpha\beta}^{\gamma}(t) & \longrightarrow &
\Gamma_{\alpha\beta}^{\gamma}(t) +
\xi_{\alpha}(t)\,\delta_{\beta}^{\gamma} +
\xi_{\beta}(t)\,\delta_{\alpha}^{\gamma}.
\end{eqnarray*}

Therefore we conclude that the neighbourhood $M_{0}$ of the point $t_{0}$
in the  moduli space comes equipped {\em canonically} with a
projective structure.

Let us now prove that for each point $y_{0}\in Y'= \cup_{t\in M_{0}}X_{t}$,
the associated submanifold $P_{y} \subseteq
\bnu\circ \bmu^{-1}(y) \subset M_{0}$ is totally geodesic relative to the
canonical projective connection in $M_{0}$.
Suppose that $y_{0}\in W_{\i}$ for some $\i$. Then $y_{0} = (w_{\i\,0},
z_{\i\,0}^{a})$ and the submanifold $P_{y_{0}}$
is given locally by the equations
$w_{\i\,0} - \phi_{\i}(z_{\i\,0}, t) = 0$,
where $t\in \bnu\circ\bmu^{-1}(y_{0})\setminus\{\mbox{\rm singular points}\}$.
Then a vector field $\left.v(t) = V^{\alpha}\p_{\alpha}\right|_{P_{y_{0}}}$ is
tangent to $P_{y_{0}}$ if and only if it satisfies
the simultaneous  equations
\begin{equation}
V^{\alpha}\,\partial_{\alpha}\phi_{\i}(z_{\i\,0},t) = 0 . \label{010}
\end{equation}
In order to prove that the submanifold $P_{y_{0}}$ for arbitrary $y_{0}\in
Y'$
is totally geodesic relative to the canonical projective connection,
we have to show
that, for any vector fields $v(t) = V^{\alpha}\p _{\alpha}$ and $w(t) =
W^{\alpha}\p _{\alpha}$
on $P_{y_{0}}$, the equation
\begin{equation}
\left( W^{\beta}\p_{\beta}V^{\alpha} +
\Gamma_{\beta\gamma}^{\alpha}\,V^{\gamma}\,W^{\beta}\right)\bmod TP_{y_{0}} =
0.
                                                                 \label{011}
\end{equation}
holds.
Since  $v(t)$ and $w(t)$
are tangent to $P_{y_{0}}\subset M$, we have
the equation
\begin{equation}
W^{\beta}(t)\,\frac{\p}{\p t^{\beta}}\left( V^{\alpha}\,
\partial_{\alpha}\phi_{\i}(z_{\i\,0},t)\right. = 0.  \label{012}
\end{equation}
$$
V^{\alpha}\,W^{\beta}\,
\frac{\partial^{2}\phi_{\i}(z_{\i\,0},t)}{\partial t^{\alpha}\partial
t^{\beta}} = V^{\alpha }W^{\beta }\Gamma_{\alpha\beta}^{\gamma}
\frac{\partial \phi_{\i}(z_{\i\,0}, t)}{\partial t^{\gamma}} .
$$
{}From the latter equation and equation (\ref{012}) it follows that
$$
\left( W^{\beta}\p_{\beta}V^{\alpha} +
\Gamma_{\beta\gamma}^{\alpha}\,V^{\gamma}\,W^{\beta}\right)
\frac{\p \phi_{\i}(z_{\i\,0},t)}{\p t^{\alpha}} = 0.
$$
By (\ref{010}) this means that $\left( W^{\beta}\p_{\beta}V^{\alpha} +
\Gamma_{\beta\gamma}^{\alpha}\,V^{\gamma}\,W^{\beta}\right)\p_{\alpha}
\in TP_{y_{0}}$, and thus  equation (\ref{011}) holds. The proof is completed.
\, $\Box$

We may have a moduli space even if the condition $H^{1}(X, N) = 0$ is not
satisfied. Given a moduli space, the proof above provides a projective
structure so we have the following global result.

\begin{corollary}\label{th 2}
Let \mbox{$\{X_t\hook Y \mid t\in M \}$} be a complete analytic family
of compact complex hypersurfaces such that $H^{1}(X_t,\f_{X_t}) = 0$ for
all $t\in M$. Then the moduli space $M$ comes
equipped canonically with a projective structure
such that, for every point $y\in Y' $, the associated submanifold $P_{y} =
\bnu\circ \bmu^{-1}(y) \subset M$ is totally geodesic.
\end{corollary}

We conclude this section with a brief geometric interpretation of geodesics
canonically induced on moduli spaces of compact complex
hypersurfaces. Any complex curve (immersed connected complex 1-manifold) in a
complex manifold $M$  has a canonical lift to a complex curve in the
projectivized tangent bundle $P_M(TM)$ --- one simply associates to each point
of the curve its tangent direction. Then a projective structure on $M$
defines a family of lifted curves in $P_M(TM)$ which foliates the
projectivized bundle holomorphically \cite{Hitchin,L}. Then, for geodesically
convex $M$, the quotient space of this foliation, $Z$, is a
$(2n-2)$-dimensional manifold, where $n=\dim M$. There is a double fibration
\begin{equation}
Z \stackrel{\tau}{\longleftarrow} P_M(TM) \stackrel{\sigma}{\lon} M
\label{f}
\end{equation}
such that, for each $z\in Z$, $\sigma\circ\tau^{-1}(z) \subset M$ is a
geodesic from the projective structure; for each $t\in M$, $\tau\circ
\sigma^{-1}(t)\subset Z$ is projective space $\Bbb C\Bbb P^{n-1}$ embedded
into $Z$ with normal bundle $T\Bbb C\Bbb P^{n-1}(-1)$ \cite{L}.

Let \mbox{$X_0 \hookrightarrow Y$} be a compact complex submanifold of
codimension 1 such that $H^{1}(X_0, N) = H^1(X_0,\f_{X_0})= 0$ and let $M$ be
a geodesically convex domain in the associated complete moduli space of
relative deformations of $X_0$ inside $Y$. The space of geodesics $Z$ can be
identified in this case with the family of intersections
$X_s\cap X_t\subset Y$, $s,t\in M$. From the explicit coordinate description
of submanifolds $X_t\subset Y$ given in the proof of Theorem\ref{th 1}
one can easily see that, for each $t\in M$,
the intersection $X_t\cap X_0$ is a divisor of the holomorphic line
bundle on $X_0$ which is a holomorphic deformation of the normal bundle $N$.
Since $H^1(X_0,\f_X)=0$, any holomorphic deformation of $N$ must be isomorphic
to $N$ \cite{Kodaira3}. Therefore, each intersection $X_t\cap X_0$ is a
divisor of the normal bunlde on $X_0$, and, by completeness
of the family $\{X_t\hook Y\mid t\in M\}$, all divisors of $N$ arise in this
way. If $t_0\in M$ is the point associated to $X_0\subset Y$ via the double
fibration (\ref{F1}), then the set of all intersections
$X_0\cap X_t$ is a projective space
$\Bbb C \Bbb P^{\dim M-1}\subset Z$ associated to $t_0$ via the double
fibration~(\ref{f}). Then a geodesic
through the point $t_0\in M$ is a family
of $X_t$ which have the same intersection with $X_0$.

\paragraph{4. Applications and examples.}
One of the immediate applications of the
theorem on projective connections is in the theory of 3-dimensional
Einstein-Weyl manifolds. Hitchin \cite{Hitchin} proved that there is a
one-to-one correspondence between local solutions of Einstein-Weyl equations in
3 dimensions and pairs $(X,Y)$, where $Y$ is a complex 2-fold and $X$ is the
projective line $\Bbb C\Bbb P^{1}$ embedded into $Y$ with normal bundle
$N\simeq \f (2)$. However the corresponding twistor techniques allowed one to
compute only part of the canonical Einstein-Weyl structure induced on the
complete moduli space $M$ of relative deformations of $X$ in $Y$, namely the
conformal structure on $M$. Although the geodesics were formally described and
the existence of a connection with special curvature was proved, no explicit
formula for the connection was obtained.  The theorem on projective connections
fills this gap and provides one with a technique which is capable of decoding
the full Einstein-Weyl structure from the holomorphic data of the embedding
$X\hookrightarrow Y$. We shall use Theorem~\ref{th 1} in some examples to
compute explicitly the canonical projective connection and then the canonical
Einstein-Weyl structure on the complete moduli space of rational curves
embedded into a 2-dimensional complex manifold with normal bundle $N\simeq \f
(2)$.

Consider a non-singular curve $X$ of bidegree $(1,n)$ in the quadric
${\Bbb C}{\Bbb P}^1 \times {\Bbb C}{\Bbb P}^1$. Then $X$ is rational and
has normal bundle ${\cal O}(2n)$ \cite{Pedersen}. The space $M$ of such
curves can be described as follows: Let $(\zeta ,\eta )$ be affine
coordinates on ${\Bbb C}{\Bbb P}^1 \times {\Bbb C}{\Bbb P}^1$ and
consider the graph of a rational function of degree $n$:
\begin{eqnarray}
\eta = \frac{P(\zeta )}{Q(\zeta )} \nonumber \\
P(\zeta )=a_n\zeta^n+a_{n-1}\zeta^{n-1}+ \cdots +a_0  \label{HP}\\
Q(\zeta )=b_n \zeta^n +b_{n-1}\zeta^{n-1}+ \cdots +b_0 \nonumber
\end{eqnarray}
The family of such $(1,n)$-curves is parameterized by ${\Bbb C}{\Bbb
P}^{2n+1}$ and the space $M$ of non-singular curves is ${\Bbb C}{\Bbb
P}^{2n+1} \backslash R$ where $R$ is the manifold of codimension 1 and
degree $2n$
given by the resultant of $P$ and $Q$. The geodesics of the projective
connection are again given by projective lines in ${\Bbb C}{\Bbb
P}^{2n+1}\backslash R$. We may of course choose to describe the induced
structure on the hypersurface given by $R=1$, and for $n=1$ this corresponds to
the standard projective structure on $SL(2,{\Bbb C})$ or on one of its real
slices ${\cal H}^3$, ${S}^3$.

In order to obtain less trivial examples we consider branched coverings.
Consider a complex curve $C$ contained in a complex surface ${S}$.
We want to construct a branched covering of a neighborhood of $C$
branched along $C$. Choose coordinates $(x_{\i} , y_{\i})$ on neighborhoods
$O_{\i}$ along $C$ such that $O_{\i} \cap C$ is given by $x_{\i} =0$.
Then on overlaps we have
$x_{\i} =x_{\j} H_{\i\j}(x_{\j} , y_{\j})$ and
$y_{\i} =K_{\i\j}(x_{\j} , y_{\j})$.
Now, we look for an $n$-fold cover branched along $C$: take patches $W_{\i}$
with coordinates $(w_{\i},z_{\i})$ and define the covering map
$(w_{\i} , z_{\i}) \to (x_{\i} , y_{\i})=(w_{\i}^n,z_{\i})$.
This is a branched cover of $O_{\i}$ branched along $O_{\i} \cap
C$. We want to identify the neighborhoods $W_{\i}$ along the curve $C$ to
obtain a surface $Y$ with a map $\pi :Y \to S$ which locally has the
form above. We get
\[
w_{\i}^n=x_{\i}=x_{\j}
H_{\i\j}(z_{\j},y_{\j})=w_{\j}^nH_{\i\j}(w_{\j}^n,z_{\j})
\]
If we make a choice of the $n$-th root and put $\stackrel{\sim
}{H}_{\i\j}=H^{\frac{1}{n}}_{\i\j}$ we get
\[
w_{\i}=w_{\j} \stackrel{\sim }{H}_{\i\j}(w_{\j}^n,z_{\j})=f_{\i\j}(w_{\j},
z_{\j}).
\]
The obstruction for this to work along the curve is the class
$\stackrel{\sim }{H}_{\i\j} \stackrel{\sim }{H}_{\j\k} \stackrel{\sim
}{H}_{\k\i} \in H^2(C,{\Bbb Z}/n)$.
We can identify this obstruction to be the self-intersection number of
$C$ modulo $n$: since $d x_{\i} =H_{\i\j}(0,y_{\j})d x_{\j}$ we see that
$H_{\i\j}(0,y_{\j})$ represents the normal bundle $N$ in $H^1(C,{\cal O}^{\ast
})$. From the long exact sequence associated with
\[
0 \to {\Bbb Z} \to {\cal O} \stackrel{\exp}{\to } {\cal O}^{\ast } \to 0
\]
we see that the degree of $N$ is equal to $log \ H_{\i\j}+ \ log \ H_{\j\k}
+ \ log H_{\k\i}$. Thus, the obstruction to obtain Y is equal to the
self-intersection of $C$, modulo $n$. Each choice of
$H^{\frac{1}{n}}_{\i\j}$ corresponds to an element in $H^1(C,{\Bbb
Z}/n)$. Unless the homology class of $C$ in $H_2(S,{\Bbb Z})$ is
divisible by $n$ this local construction along the curve cannot be
extended to work globally on $S$ \cite{Atiyah}.

Now, let us return to the case where $C$ is a $(1,n)$--curve in ${\Bbb
C}{\Bbb P}^1 \times {\Bbb C}{\Bbb P}^1$. In this case $C \cong {\Bbb
C}{\Bbb P}^1$, so there is a unique $n$--fold covering $Y$ branched
along $C$ which we cannot extend to all of ${\Bbb C}{\Bbb P}^1 \times
{\Bbb C}{\Bbb P}^1$. The branch locus $X \subseteq Y$ is a copy of $C$
but $deg \ N_X=\frac{1}{n} \ deg \ N_C =2$
so we may describe an Einstein-Weyl structure on the moduli space of
curves in $Y$ \cite{Hitchin} and contrary to earlier attempts we are now able
to get
the connection $\Gamma^{\gamma}_{\alpha \beta }$ explicitly. Let us
concentrate on $(1,2)$ curves and let $C$ be the curve
$\eta = \zeta^2$.
The projection $\pi $ maps the curves in $Y$ onto those $(1,2)$--curves
which meet $C$ in two points to second order. These curves may be given
as in (\ref{HP}) with
\[
P(\zeta )=\zeta^2-2t_0t_1\zeta-t_0^2
\]
\[
Q(\zeta )=t^2_2\zeta^2+2t_1t_2\zeta +1+2t_0t_2+t_1^2
\]
(see \cite{Pedersen}). In order to describe the lifted curves we introduce
the coordinates
\[
\begin{array}{ll}
x_1=\eta -\zeta^2 & x_2 = \stackrel{\sim }{\eta }-\stackrel{\sim }{\zeta }^2 \\
y_1=\zeta & y_2=\stackrel{\sim }{\zeta }
\end{array}
\]
where $(\stackrel{\sim }{\zeta },\stackrel{\sim }{\eta
})=(\frac{1}{ \zeta },\frac{1}{\eta })$. Then $C$ is given by $x_{\i}=0$.
Making the coordinate transformation
\begin{eqnarray*}
(x_{1}, y_{1}) & \lon & (w, z) = (\sqrt{x_{1}}, y_{1}) \\
(x_{2}, y_{2}) & \lon & (\hat{w}, \hat{z}) = (\sqrt{x_{1}}, y_{1})
\end{eqnarray*}
we arrive at a covering of $Y$ by two coordinate charts $W$ and $\hat{W}$ which
is exactly of the type used in the proof of Theorem~1 and has
the transition functions
$$
\hat{w} = f(w,z), \ \ \ \  \hat{z} = g(z),
$$
given by
\begin{eqnarray*}
f(w,z) & = & \frac{w}{z\, \sqrt{w^{2} + z^{2}}}\\
g(z) & = &  z^{-1}.
\end{eqnarray*}
The complete maximal family of relative deformations of $C$ is
described in this chart
by the equations (in the notation of the proof of Theorem~\ref{th 1})
$w = \phi(z,t)$ and $\hat{w} = \hat{\phi}(\hat{z}, t)$,
with
$$
\phi(z,t) = i\,R(z)\,Q(z)^{-1/2}, \ \ \ \ \hat{\phi}(z,t) =
i\,R(z)\,P(z)^{-1/2},
$$
where $R(z) = t_{2}\,z^{2} + t_{1}\,z + t_{0}$.
Note that a useful identity $P = z^{2}Q - R^{2}$ holds \cite{Pedersen}.

Now we have all the data to apply the machinery developed in the proof
of Theorem~\ref{th 1}. Following that scenario one finds that the canonical
projective structure on $M$ can be represented by the following torsion-free
affine connection
$$
\begin{array}{rclrcl}
\vspace{4 mm}
\Gamma_{01}^{0} & = & t_{1}\,(1 + 3\,t_{0}\,t_{2})\,(2\,\triangle)^{-1}, &
\Gamma_{01}^{1} & = & t_{2}\,(2 + t_{1}^{2} +
2\,t_{0}\,t_{2})(2\,\triangle)^{-1}, \\
\vspace{4 mm}
\Gamma_{00}^{0} & = & t_{2}\,(1 + t_{0}\,t_{2})\,\triangle^{-1},&
\Gamma_{00}^{1} & = & - t_{1}\,t_{2}^{2}\,\triangle^{-1} , \\
\vspace{4 mm}
\Gamma_{02}^{0} & = & t_{0}\,( 1 + t_{0}\,t_{2} +
t_{1}^{2})(2\,\triangle)^{-1}, &
\Gamma_{02}^{1} & = & - t_{1}\,(1 + \,t_{1}^{2})(2\,\triangle)^{-1},\\
\vspace{4 mm}
\Gamma_{11}^{0} & = & - t_{0}\,(1 + t_{0}\,t_{2})\,\triangle^{-1}, &
\Gamma_{11}^{1} & = & t_{0}\,t_{1}\,t_{2}\,\triangle^{-1}\\
\Gamma^{0}_{12} & = & - t_{0}^{2}\,t_{1}\,(2\,\triangle)^{-1}, &
\Gamma^{1}_{12} & = & - t_{0}\,( 1 + t_{0}\,t_{2} + t_{1}^{2})
(2\,\triangle)^{-1}
\end{array}
$$
and all other Christoffel symbols being zero. Here
$$
\triangle = (1 + t_{0}\,t_{2})^{2} + t_{1}^{2}(1 + 2\,t_{0}\,t_{2}).
$$
Note that $\triangle^2=R$ where $R$ is the resultant of the polynomials
in (\ref{HP}).

The conformal structure $[g]$ on $M$ is given by the condition for the
curves to meet to second order. Thus we may choose
the following metric in the conformal structure \cite{Pedersen}
\begin{eqnarray}
g &=& t^2_1t^2_2dt^2_0+(1+t_0t_2)^2dt^2_1+4t^2_0(1+t^2_1)dt^2_2 +
2t_1t_2(1+t_0t_2)dt_0dt_1 \nonumber \\
&-&4(1+t^2_1)(1+t_0t_2)dt_0dt_2  - 4t_0^2t_1t_2dt_1dt_2
\end{eqnarray}
Since our connection $\nabla $ is projectively equivalent to
the Weyl connection $D$ it satisfies
\[
(\nabla g)_{\alpha \beta \gamma } = a_{\alpha }g_{\beta \gamma
}+b_{\beta }g_{\alpha \gamma }+b_{\gamma }g_{\alpha \beta }
\]
for some 1-forms $a=\sum_{i=0}^{2}a_{\alpha}\,dt^{\alpha}$ and $b =
\sum_{i=0}^{2}b_{\alpha}\,dt^{\alpha}$.
We may solve these equations and present the Weyl connection $D$ in
terms of the Levi-Civita connection $\nabla^g$ and the 1-form
$\omega = a - 2b =\sum_{\alpha }\omega_{\alpha }dt_{\alpha }$,
\[
D=\nabla^g+\frac{1}{2} \omega^{\# }g-\omega \odot I
\]
see \cite{Pedersen-Tod}. We get

$$
\begin{array}{rclrcl}
\vspace{4 mm}
a_{0} & = & 3\,t_{1}^{2}\,t_{2}\,(2\,\triangle)^{-1}, &
a_{1} & = & - 3\,t_{1}\,(1 + t_{0}\,t_{2})(4\,\triangle)^{-1}, \\
\vspace{4 mm}
a_{2} & = & - 3\,t_{0}\,(1 + t_{0}\,t_{2} + t_{1}^{2})(2\,\triangle)^{-1},
& & & \\
\vspace{4 mm}
b_{0} & = & - 3\,t_{1}^{2}\,t_{2}\,(4\,\triangle)^{-1}, &
b_{1} & = & - 3\,t_{1}\,(1 + t_{0}\,t_{1})(4\,\triangle)^{-1}, \\
b_{2} & = & - 3\,t_{0}\,(1 + t_{0}\,t_{2} + t_{1}^{2})(2\,\triangle)^{-1}.
&&&
\end{array}
$$
Thus using only the methods of the relative deformation theory of compact
hypersurfaces we computed the full Einstein-Weyl structure on the
moduli space.

Suppose we blow up a point $s$ on the quadric and take a $(1,n)$-curve
passing through the point. Then in the blown up surface the curve will
have self-intersection number $2n-1$ and this corresponds to considering
all the $(1,n)$-curves passing through $s$. We may combine this with
the branched covering construction. In \cite{Pedersen-Tod} we considered
the Einstein-Weyl structure associated to the $(1,3)$-curves: first we
considered the 2-fold branched cover which reduced the degree of the
normal bundle from 6 to 3 and then we blew up a point on the branch
locus to get self-intersection equal to 2. Again we may compute the Weyl
connection or compute the connection associated to any combination of
blow up and branched cover. This will give non trivial examples with
normal bundle ${\cal O}(n)$ for any $n$.
\vspace{2 mm}

{\em Acknowledgments}. It is a pleasure to thank Paul Tod for many valuable
discussions and comments. Thanks are also due to Stephen Huggett, Yat Sun Poon
and the anonymous referees for helpful remarks. One of the authors (SM) is
grateful to the Department of Mathematics and Computer Science of Odense
University for hospitality and financial support.

\pagebreak

\end{document}